\documentclass[11pt,ams]{article}
\usepackage{amssymb, amscd, amsmath, amsthm}

\newtheorem{thm}{Theorem}
\newtheorem{question}[thm]{Question}
\newtheorem{lem}[thm]{Lemma}
\newtheorem{prop}[thm]{Proposition}
\newtheorem{cor}[thm]{Corollary}
\allowdisplaybreaks

\def\Sym{\text{\rm Sym}\,}
\def\ol{\overline }
\begin{document}

\title{Product decompositions of quasirandom groups and a Jordan
type theorem}

\author{N. Nikolov and L. Pyber}

\date{}

\maketitle

\begin{abstract}
We first note that a result of Gowers on product-free sets in
groups has an unexpected consequence:
If $k$ is the minimal degree of a representation of the finite
group~$G$, then for every subset $B$ of $G$ with $|B| > |G| /
k^{\frac13}$ we have $B^3 = G$.

We use this to obtain improved versions of recent deep theorems
of Helfgott and of Shalev concerning product decompositions of
finite simple groups, with much simpler proofs.

On the other hand, we prove a version of Jordan's theorem which
implies that if $k \geq 2$, then $G$ has a proper subgroup of
index at most $c_0 k^2$ for some constant $c_0$, hence a product-free subset of size at least 
$|G| / ck$.
This answers a question of Gowers.
\end{abstract}

\setcounter{section}{-1}
\section{Introduction}
\label{sec:0}

Sum-free subsets of abelian groups have been much investigated
in the past 40 years. Very recently Green and Ruzsa \cite{GR} have
determined the maximal size of a sum-free subset of any finite abelian group.

A subset $X$ of a not necessarily abelian  group $G$ is called
\emph{product-free} if there are no solutions to $xy = z$ with
$x, y, z \in X$.
The maximal size $\alpha(G)$ of a product-free subset of a
finite group $G$ has been considered by Babai and S\'os \cite{BS}.
For example, they proved that for any soluble group $G$ of order
$n$ we have $\alpha(G) \geq \frac{2n}{7}$ and asked whether a
similar linear bound holds for arbitrary finite groups.

Note that any non-trivial coset of a subgroup is product-free.
In fact, in 1950 (according to \cite[p.~246]{PS}) Moser
conjectured that the largest product-free subsets of alternating
groups are cosets of maximal subgroups.

Kedlaya \cite{Ke} disproved this, by showing that if a subgroup
$H$ has index $m$ in a group~$G$, then one can actually find a
set of size $cm^{1/2}|H|$ that is product-free. Here and below $c$ denotes an absolute constant.

Combining this estimate with the classification of finite simple
groups (CFSG) Kedlaya showed that for every finite group~$G$ we
have $\alpha(G) \geq cn^{\frac{11}{14}}$.
He asked whether for every $\varepsilon > 0$ one can obtain a
bound of $c(\varepsilon)n^{1 - \varepsilon}$.

A negative answer to the above question was obtained very
recently by Gowers \cite{Gow}.
He showed that for sufficiently large $q$ the group $\Gamma =
PSL(2,q)$ has no product-free subsets of size $c|\Gamma|^{\frac89}$
His proof depends on the fact, proved by Frobenius, that every
non-trivial representation of $PSL(2, q)$ has degree at least
$(q - 1)/2$.

Gowers went on to consider combinatorial properties of finite
groups $G$ such that every non-trivial representation of $G$ has
degree at least~$k$.
He calls such groups \emph{quasi\-random}, since this property
turned out to be equivalent to several other properties, some of
which state that certain associated graphs are quasi\-random (see
\cite{Gow} for a detailed discussion of quasi\-random graphs).

Gowers proved \cite[p.~22]{Gow} the following general result on product-free sets
of quasi\-random groups.

\setcounter{thm}{-1}

\begin{prop}
\label{prop:0}
Let $G$ be a group of order $n$, such that the minimal degree of
a nontrivial representation is~$k$.
If $A, B, C$ are three subsets of~$G$ such that $|A|\, |B|\, |C|
> \frac{n^3}{k}$, then there is a triple $(a, b, c) \in A \times
B \times C$ such that $ab = c$. \hfill $\square$
\end{prop}

The starting point of the present paper is the following
surprising consequence.

\begin{cor}
\label{cor:1}
Let $G$ be a group of order $n$, such that the minimal degree of
a representation is~$k$.
If $A, B, C$ are three subsets of $G$ such that $|A|\, |B|\, |C|
> \frac{n^3}{k}$, then we have $A \cdot B \cdot C = G$.
In particular, if, say, $|B| > \frac{n}{k^{1/3}}$, then we have
$B^3 = G$.
\end{cor}

\begin{proof}
Consider the set $G \setminus AB$.
By Proposition~\ref{prop:0} the size of this set is strictly
less than $|C|$, i.e., we have $|AB| + |C| > |G|$.
It follows that for any $g \in G$ the intersection of the sets
$AB$ and $gC^{-1}$ is non-empty, which implies $g \in ABC$.
\end{proof}

Applying Corollary~\ref{cor:1} to the sets $A, B, C^{-1}$ we see
that Proposition~\ref{prop:0} and Corollary~\ref{cor:1} are in
fact equivalent.

Corollary~\ref{cor:1} apart from its intrinsic interest, seems
to be an extremely useful tool.
Recently a number of deep theorems have been obtained concerning
product decompositions of simple groups.
Corollary~\ref{cor:1} can be used to give short and relatively
elementary proofs, while improving the results.

It is particularly useful in the case of simple groups of Lie type.
For these groups rather strong lower bounds on the minimal
degree of a representation are known.

For a (possibly twisted) Lie type $L$, not $^2B_2$, $^2G_2$,
$^2F_4$ define the rank $r = r(L)$ to be the untwisted Lie rank
of $L$ (that is, the rank of the ambient simple algebraic group)
and for $L$ of type $^2B_2$, $^2G_2$, $^2F_4$ define $r(L) = 1,
1, 2$ respectively.

It follows from \cite{LS} that there is a constant $c$ such that
for any simple group $L$ of Lie type of rank $r$ over $\mathbb
F_q$ we have $k \geq cq^r$ for the minimal degree $k$ of a
representation of~$L$.

For $L = PSL(n, q)$ we obtain the following.

\begin{prop}
\label{prop:2}
Let $B$ be a subset of $L$ of size at least
${2|L|}/{q^{\frac{n - 1}{3}}}$.
Then we have $B^3 = L$.
\hfill $\square$
\end{prop}

A similar result in the case of $\Gamma = PSL(2, p)$,
$p$ prime, plays an important role in the proof of a recent
breakthrough result of Helfgott concerning the diameter of
Cayley graphs of~$\Gamma$.

Helfgott \cite{He} showed that for every set of generators $X$
of $\Gamma$ every element of $\Gamma$ can be expressed as a
product of at most $O((\log p)^c)$ elements of $X \cup X^{-1}$.

Proposition~\ref{prop:2} improves (the easier part of) his Key
Proposition \cite[p.~2]{He} even for $\Gamma= PSL(2,p)$ and implies that
in fact every element of $\Gamma$ can be expressed as a product
of at most $O((\log p)^c)$ elements of $X$ itself (this improvement 
also follows from the results in \cite{Babai}).

It is an open problem whether Helfgott's result extends to
$PSL(n, p)$.
Proposition~\ref{prop:2} may be useful in obtaining a positive answer.

As another interesting application we prove the following Waring
type theorem. For a group word $w = w(x_1, \dots, x_d)$ let $w(G)$
denote the set of values of $w$ in~$G$.

\begin{thm}
\label{thm:3}
Let $k \geq 1$ and $\ol w = \{w_1, \dots, w_k\}$ be a set of
non-trivial group words.
Let $L$ be a finite simple group of Lie type of rank $r$ over the field $F_q$ and set
$\ol w(L) = w_1(L) \cap \dots \cap w_k(L)$. Let $W$ be any subset of $\ol w(L)$ such that
$|W| \geq | \ol w(L) |/q^{r/13}$.

There exists a positive integer $N$ depending only on $\ol w$
such that if $|L| > N$, then we have
\[
W^3 = L . \eqno{\square}
\]
\end{thm}

As the main result of a difficult paper Shalev \cite{Sh} has
obtained the same result in the case $k = 1$ and $W=\ol w(L)$ (allowing $L$ to be also an alternating group). 
\footnote{We remark that for alternating groups and simple groups of Lie type of bounded rank it was shown later in
\cite{LSh1} that in fact one has $\ol w (L)^2 = L$
if $L$ is large enough.}

Combining the methods of that paper with \cite{LSh1} one can prove it for
$W=\ol w(L)$ and $k$ arbitrary. An advantage of our sparse version is that one can
impose further restrictions on $W$. For example one can require that no two
elements of $W$ are inverses of one another or images of one another under
Frobenius automorphisms. In addition, (using the first part of Corollary 
\ref{cor:1}) it follows that every element $g \in L$ is a product $g=h_1 h_2 h_3 $ of \emph{distinct} 
$h_i \in \ol w (L)$. It will be interesting to see if $h_i$ can be taken to be pairwise noncommuting elements from $\ol w (L)$, or such that $\langle h_1,h_2,h_3 \rangle =G$; this doesn't seem to follow immediately from Corollary \ref{cor:1}. 
\medskip

Shalev's proof in \cite{Sh} relies on a
whole array of deep results on the character theory of simple
groups, developed to estimate the diameters of Cayley graphs of
simple groups with respect to conjugacy classes.

Our proof of Theorem \ref{thm:3} is relatively short compared to \cite{Sh} 
and uses an auxiliary result from \cite{LSh1}, see Proposition \ref{prop:1.2} below.
This says roughly that for simple groups of Lie type not of type
$A_r$ or $^2A_r$ the sets $\ol w(L)$ are ``very large''.

For groups of type $A_r$ and $^2A_r$ we provide somewhat weaker
estimates for $|\ol w(L)|$ which still make
Corollary~\ref{cor:1} applicable.

It would be most useful to obtain analogues of
Corollary~\ref{cor:1} for smaller sets~$B$.
The following results indicate how far one can go in this
direction.

\begin{thm}
\label{thm:4}
Let $G$ be a finite linear group of degree $k$ over the complex field.
Then $G$ has a permutation representation of degree at most
$ c_0 k^2$ with abelian kernel, where $c_0<10^{10}$ is an absolute constant.\hfill $\square$
\end{thm}

The proof of this result relies on the Classification of the finite simple groups (CFSG).
As an immediate consequence we obtain the following.

\begin{cor}
\label{cor:5}
Let $G$ be a finite group such that $G$ has an irreducible
representation of degree $k \geq 2$.
Then $G$ has a proper subgroup $H$ of index at most~$c_0k^2$.\hfill
$\square$
\end{cor}

In particular $H$ is a subset of size at least $\frac{n}{c_0k^2}$ in
$G$ which does not even generate~$G$.

As a ``partial converse'' to Proposition~\ref{prop:0} Gowers
\cite{Gow} proved that if a group $G$ contains no large
product-free subsets, then it is quasi\-random.
More precisely, he gave an elementary argument showing that if
the minimal degree of a representation of $G$ is $k$, then $G$
has a product-free subset of size at least $\frac{n}{c^k}$ for
some absolute constant $c > 1$.
Gowers asked whether this can be improved to $\frac{n}{k^c}$ (for $k
\geq 2$).
Applying Kedlaya's result to $H$ as above we see that $G$ has
a product-free subset of size at least $\frac{n}{ck}$ for some constant $c$, i.e.\ we
obtain a positive answer to his question. \medskip

Finally, for completeness in the last section we present a simplified version of Gowers' proof 
of Proposition \ref{prop:0} in the special case when one 
of the sets $A$, $B$ or $C$ is symmetric. This case is enough for most of our applications above.

\section{Waring type theorems}
\label{sec:1}

The main result of \cite{Sh} is the following

\renewcommand{\thethm}{\thesection.\arabic{thm}}

\setcounter{thm}{0}

\begin{thm}
\label{thm:1.1}
Let $w \neq 1$ be a group word.
Then there exists a positive integer $N = N(w)$ such that for
every nonabelian finite simple group $G$ with $|G| \geq N$ we have
\[
w(G)^3 = G.
\eqno{\square}
\]
\end{thm}

For example, each $g \in G$ can be expressed as a product of
three $k$-th powers.

The proof in \cite{Sh} relies on algebraic geometry via \cite{La}, the
Deligne--Lusztig theory of characters of Chevalley groups and on
a recent work on character theoretic zeta functions \cite{LSh2}.

In this section we indicate a proof of Theorem \ref{thm:3}. It is a
generalization of the above theorem for groups of Lie type which
does not use such difficult character theoretic tools.
We rely instead on some very recent results of Larsen
and Shalev. In \cite{LSh1} they give a short proof of the following.

\begin{prop}
\label{prop:1.2} Let $k \geq 1$ and $\ol w = \{w_1, \dots, w_k\}$
be a set of words. Let $L$ be a finite simple group of Lie type of
rank $r$, which is not of type $A_r$ or $^2A_r$. There is an
absolute constant $c$ and an integer $N=N(\ol w)$ such that if $|L| \geq N$ then we
have
\[
|\ol w (L)| \geq c|L|/r.
\eqno{\square}
\]
\end{prop}

The proof of Proposition~\ref{prop:1.2} relies on algebraic
geometry (replacing \cite{La}) and group theoretic arguments.
It also applies to groups of type $A_r$ and $^2A_r$ for $r$
bounded and in fact to their covering groups $SL(r + 1, q)$,
$SU(r + 1, q)$.

As noted in the introduction, if $L$ is a simple group of Lie
type over $\mathbb F_q$, then we have $k \geq c q^r$ for the
minimal degree $k$ of a representation of~$L$.
Hence Proposition~\ref{prop:1.2} and Corollary~\ref{cor:1}
immediately imply Theorem~\ref{thm:3} if $L$ is not of type
$A_r$ or $^2A_r$ with $r \geq 1000$, say.

For the groups $SL(n, q)$ and $SU(n,q)$  the minimal degree of a
representation is at least $\frac{q^{n - 1} - 1}{2}$, see \cite{LS}.
To complete the proof of Theorem~\ref{thm:3} we have to show
that for these groups we have
\[
\frac{|\ol w(G)|}{q^{\frac{n-1}{13}}} > \frac{|G|}{\left(\frac{q^{n - 1} - 1}{2}\right)^{\frac13}}
\]
if $n$ is large enough. This is achieved below where we prove that
$|\ol w (G)|/|G|  \geq \gamma n^{-3}q^{-50 -n/4}$ for some 
constant $\gamma>0$ depending only on $\ol w$.
\bigskip

Let $G = L(q)$ be a quasisimple group of Lie rank $r$ defined over
$\mathbb F_q$. A regular semisimple (r.s.) element $g$ of $G$
is one which has distinct eigenvalues (possibly in a field
extension of $\mathbb F_q$). Therefore $C_G(g)$ is a torus of
$G$ and so $|C_G(g)|=(1+o(1))q^r$.  It is well-known that the
set of regular semisimple elements of a semisimple connected
algebraic group $G$ has complement of strictly smaller dimension
than $\dim G$ (see \cite{GL2}).
Using this we have

\begin{lem}
\label{reg}
Given the type $L$ there is a constant $C$ depending on $L$ such
that the cardinality of the r.s.\ elements of $G$ is at least
$(1-C/q)|G|$.
\end{lem}

The following proposition is an important ingredient for our proof
of Theorem~\ref{thm:3}.

\begin{prop}
\label{sl4}
Given  $\ol w = \{ w_1, \dots, w_k\} $ there is a
constant $c>0$ depending only on $\ol w$
such that if
$G=SL(4n,q)$, then $| \ol w(G)|> c |G|/(n^3 q^{n-1})$. In fact,
$\ol w(G)$
contains at least $cq^{3n}/n^3$ conjugacy classes of regular
semisimple elements.
\end{prop}

\begin{proof}[Proof of Proposition \ref{sl4}]
The main idea of the proof is a generalization of some arguments
in Section 2 of \cite{LSh1} from the case of $SL(2,q^n)$ to
$SL(4,q^n)$:

Consider the inclusion $i: H=SL(4,q^n) \rightarrow SL(4n,q)$ and
let $g \in H$ be a semisimple element whose eigenvalues form the
multiset $\{a_1, \ldots , a_4\}=A$. Let $F$ be the automorphism
$x \mapsto x^{q}$ of $\overline{ \mathbb F}_q$. Since the
eigenvalues $a_i$ are roots of the characteristic polynomial of
$g$ with coefficients in $\mathbb F_{q^n}$ it follows that
$A^{F^n}$ is a permutation of $A$.

The eigenvalues of $i(g)$ form the multiset $A, A^F, \ldots
A^{F^{n-1}}$.

\begin{lem}\label{ss}
 There are at least $(1-O(nq^{-n/2}))
|H|$ elements of $g \in H$ such that $i(g)$ is regular
semisimple in $G$.
\end{lem}

\begin{proof}
By Lemma \ref{reg} the number of elements of $H$
which are not regular semisimple in $H$ is $O(q^{-n})|H|$.
So it is enough to consider only regular semisimple elements $g\in H$.

Now, suppose that $g \in H$ is regular semisimple but $i(g) \in
G$ is not regular. The eigenvalue multiset of $g \in H$ then
consists of 4 distinct elements and has the form
\[
A=\bigcup_{i=1}^s \{\alpha_i, \alpha_i^{F^{n}},\cdots
,\alpha_i^{F^{(k_i-1)n}} \},
\]
where $4=k_1+\cdots k_s$ is a partition of $4$, and $\alpha_i$
is a generator for the finite field $\mathbb F_{q^{k_i n}}$ over 
$\mathbb F_{q^n}$.  In
addition, the product of all elements of $A$ should be $1$.
\medskip

Since $i(q)$ is not regular there are two elements $\alpha$ and
$\beta$ from $A$ such that $\alpha= \beta^{F^j}$ where $0<j<n$.
Now $\alpha$ and $\beta$ are either from distinct orbits of
$F^n$ in $A$, or from the same orbit.  In the first case without
loss of generality we may assume that $\alpha =\alpha_1$ and
$\beta=\alpha_2=\alpha_1^j$. Then $k_1=k_2$ and the eigenvalue
multiset $A$ of $g$ is then determined by $j$ and the $(s-1)$
eigenvalues $\alpha_i \in \mathbb{F}_{q^{k_in}}$ for $i \not = 2$.
Simple
calculation shows that the number of possibilities for $A$ (i.e.\ for the
conjugacy class of $g$ in $H=SL(4,q^n)$) is at most
$O(nq^{2n})$.

In the second case, if we assume that $\alpha= \alpha_1$ there
is a proper divisor $d$ of $k_in$ such that
$\alpha_1=\alpha_1^{F^d}$, i.e. $\alpha_1 \in \mathbb F_{q^{d}}$. Given
$d$, counting the possibilities for $\alpha_1, \ldots \alpha_s$
under the restriction $\alpha_1 \in \mathbb F_{q^{d}}$ we see that there
are $O(q^{3n-k_in+d})=O(q^{3n-n/2})$ possible choices for $A$.
So in this case, as $d$ ranges over all proper divisors of $n$
there are altogether at most $O(nq^{5n/2})$ conjugacy classes of
such~$g$.

Combining both cases for all partitions $(k_i)$ of $4$ we see that the number of conjugacy classes
of regular semisimple elements $g \in H$ such that $i(g)$ is not
regular in $SL(4n,q)$ is $O(nq^{3n-n/2})$. This gives the
conclusion of the lemma.
\end{proof}

Now by Proposition 8.2 of \cite{LSh1} there is a constant
$c_0>0$ such that $|\ol w(H)|>c_0|H|$.
Together with Lemma \ref{ss}
this implies that $\ol w(H)$ contains at least $c_1q^{3n}$ conjugacy
classes of r.s.\ elements $g$ such that $i(g)$ is also regular.
Here the constant $c_1$ can be taken to be any number in
$(0,c_0)$ provided that $n$ is sufficiently large.

\begin{lem}
\label{glue}
Suppose $g^H$ is a conjugacy class of
r.s.\ elements of $H$. There are at most $O(n^3)$ distinct
conjugacy classes $h^H$ of $H$ such that $i(g)^G=i(h)^G$.
\end{lem}

\begin{proof}
Note that the conjugacy classes of semisimple
elements of $SL(n,q)$ are uniquely determined by the multisets of
their eigenvalues. Given $g$ with eigenvalue multiset $A$ as
above, suppose that $h$ has an eigenvalue multiset $B=\{b_1,
\ldots , b_4 \}$ and $i(g)^G=i(h)^G$, i.e.,
$\cup_{j=1}^{n-1} A^{F^j}=\cup_{j=1}^{n-1} B^{F^j}$. This
implies that every one of $b_1,b_2, b_3$ can take one of $4n$
given values. Having chosen these three, the fourth one, $b_4$
is then uniquely determined by $\det h =1$.
\end{proof}

With the above lemma we obtain that $i(\ol w(H))$ contains at least
$c q^{3n}/n^3$ distinct conjugacy classes of r.s.\ elements of $G$
(where $c$ depends only on $\ol w$) and so
\[ 
|\ol w(G)| \geq | i( \ol w(H))^G| \geq \frac{cq^{3n}}{n^3} \cdot
\frac{|G|}{q^{4n-1}}=\frac{c|G|}{n^3 q^{n-1}},\]
proving Proposition~\ref{sl4}.
\end{proof}

For general ${SL}(n,q)$ we have

\begin{prop}\label{sln}
 Given $\ol w$ there is a
constant $c'>0$ depending only on $\ol w$ and such that if
$G={SL}(n,q)$ then
\[
|\ol w(G)|>\frac{c'|G|}{n^3
q^{24+n/4}}.
\]
\end{prop}

\begin{proof}
Let $m=[n/4]$ Consider the embedding $j:
SL(4m,q) \rightarrow SL(n,q)$ in the top left corner.
By Proposition \ref{sl4} $\ol w(H)$ contains $cq^{3m}/m^3$ conjugacy
classes of r.s.\ semisimple elements $g$.
Observe that this means that $4m$ of the eigenvalues of $j(g)\in
G$ are distinct, and the rest are equal to 1. This easily gives
that $|C_G(j(g))| =O(q^{25+4m-1})$ and hence $|\ol w(G)| >
c'|G|/(n^3q^{m+24})$ for the appropriate $c'$.
\end{proof}

\subsubsection*{The unitary group}

When $L={SU}(d,q)$ the result we use is

\begin{prop}
\label{prop:1.8}
There is a constant $e>0$ such that
\[
|\ol w(L)|>\frac{e|L|}{d^3q^{49+d/4}}.
\]
\end{prop}

\begin{proof}
The argument is similar to the case of $SL(n,q)$
above. Let $n \in \mathbb{N}$ be odd. Consider the embedding $i$
of $H=SU(4,q^n)$ inside $G=SU(4n,q)$ defined in the following
way: Let $V$ be a 4-dimensional vector space over $\mathbb F_{q^{2n}}$
equipped with a nondegenerate hermitian form $v: V \times V
\rightarrow \mathbb{F}_{q^{2n}}$. Consider the form $v'= t \circ
v$ where
$t=\mathrm{Tr}_{\mathbb{F}_{q^{2n}}/\mathbb{F}_{q^{2}}}$ is the
trace map onto $\mathbb{F}_{q^{2}}$. Since $n$ is odd the above
map is still hermitian (and nondegenerate). For $g \in H$ $i(g)$
is the same transformation $g$ of $V$ considered as a vector
space over $\mathbb{F}_{q^2}$ with the form $v'$.

Let $F$ be the automorphism $x \mapsto x^q$. If $g \in H$ is a
semisimple element its multiset of four eigenvalues $A$ satisfies 
$A=A^{-F^n}$. The element 
$i(g) \in G$ has eigenvalue multiset $A,A^{F^2}, \ldots
,A^{F^{2(n-1)}}$ and from then on the argument is very much the
same: First we prove just as in Lemma \ref{ss} that the number
of elements $g \in H$ such that $i(g)$ is not r.s.\ is
$O(nq^{-n/2}|H|)$. Then, just as in Proposition \ref{sl4} it
follows that there is a constant $e_0$ such that at least
$e_0q^{3n}$ conjugacy classes of $H$ consist of elements $g \in
\ol w(H)$ such that $i(g)$ is regular semisimple. Next, at most
$O(n^3)$ such conjugacy classes of $H$ become conjugate in $G$
(because for $i(g)^G=i(h)^G$ we need three eigenvalues 
of $h$ to be from $A^{\pm F^{2j}}, \ j=0,1, \ldots ,n-1$
and then they determine the last eigenvalue of $h$). It
follows that $\ol w(G)$ contains at least $e_1q^{3n}/n^2$ conjugacy
classes of r.s.\ elements. Finally, considering the embedding of
$G={SU}(4n,q)$ in $L={SU}(d,q)$ as a subspace
group (where $n=2[(d-4)/8]+1$), we see that $\ol w(L)$ contains at
least $e_1q^{3n}/n^2$ conjugacy classes of elements $g$ with
distinct eigenvalues on a nondegenerate $4n$ dimensional
subspace $U$. Any $h \in C_L(g)$ stabilizes $U$ and $U^\perp$
and hence $|C_L(g)|=O(q^{4n+49})$. (Note that $\dim U^\perp
=d-4n\leq 7$.)
\end{proof}

This completes the proof of Theorem~\ref{thm:3}.

Incidentally we observe the following consequence of Propositions \ref{prop:1.2}, \ref{sln} and \ref{prop:1.8}.
\begin{thm}\label{gen} Given a set of words $\ol w$ and a simple group $L$ of Lie type, two random elements from $\ol w (L)$ generate $L$ with probability tending to 1 as $|L| \rightarrow \infty$. In particular a random pair of squares generates $L$ with probability 1.
\end{thm}

This follows easily from the above propositions and a result of Liebeck and Shalev in \cite{LSh3} that the set of pairs $(a,b) \in L \times L$ which don't generate $L$ has size at most $c |L|^2/P(L)<c |L|^2q^{-r}$.
 
\section{Bounds for linear groups}
\label{sec:2}

\setcounter{thm}{0}
By a classical result of Jordan a finite linear group of degree
$k$ has an abelian normal subgroup $A$ of index $j(k)$ for some
function~$j(k)$.

Essentially the best elementary estimate for $j(k)$ is due to
Blichfeldt:
\[
j(k) \leq k! \ 6^{(k - 1)\pi(k + 1) + 1}
\]
where $\pi(k + 1)$ denotes the number of primes $\leq k + 1$
(see \cite{Do}).

Better bounds can be obtained using CFSG.
Building on an unpublished work of Weisfeiler, Collins \cite{Co}
has recently shown that for $k \geq 71$ we have $j(k) = (k +
1)!$ (see \cite{GL1} for the history of this result).

It is clear that $G/A$ has an embedding into $\Sym(j(k))$.
Hence Theorem~\ref{thm:4} may be considered as a different type
of quantitative version of Jordan's theorem.

For the proof we need various auxiliary results.

A $p$-group $P$ is said to be of \emph{symplectic type} if it
has no noncyclic characteristic abelian subgroups.
The structure of these groups is well understood and is closely
related to that of extraspecial groups.

In particular if $P$ has exponent $p$ ($p$ odd), then $P$ itself
is extraspecial, and if it has exponent $4$, then $P$ is either
extraspecial or the central product of an extraspecial group
and~$\mathbb Z_4$.

We will use the following.

\begin{prop}
\label{prop:2.1}
Let $P$ be a $p$-group of symplectic type and set $C =
C_{\text{\rm Aut}(P)} (Z(P))$.

\begin{itemize}
\item[\rm (i)] If $P$ has exponent $p$ ($p$ odd) and order $p^{2m +
1}$, then $C$ can be embedded in $\Sym(p^{2m})$.

\item[\rm (ii)] If $P$ has exponent $4$ and order $2^{2m + 2}$, then
$C$ can be embedded in
$\Sym (2^{2m+2}-4)$.

\item[\rm (iii)] If $P$ has exponent $4$ and order $2^{2m + 1}$,
then $C$ can be embedded in $\Sym ( 2^{2m+1}-2)$.

\end{itemize}
Hence in all cases $C$ can be embedded in $\Sym(4p^{2m})$.
\end{prop}

\begin{proof}
The structure of the above groups $P$ and $C$ is described in
\cite[Table 4.6A]{KL}.
In all cases $\text{\rm Inn}(P)$ is an elementary abelian
minimal normal subgroup of order $p^{2m}$ in~$C$.
Each element of ${C}/{\text{\rm Inn}(P)}$ acts as a
nontrivial linear transformation, hence $\text{\rm Inn}(P)$ is
the unique minimal normal subgroup of~$C$.

In case (i) we have ${C}/{\text{\rm Inn}(P)} \cong Sp(2m,
p)$ and moreover by \cite{Gr} the extension splits.
Hence $C$ has a corefree subgroup of index $p^{2m}$ which implies~(i).

In case (ii) $P$ may be expressed as the central product of
$\mathbb Z_4$ and $m$ copies of the dihedral group $D_8$, $P =
\mathbb Z_4 \circ D_8 \circ \dots \circ D_8$.

In case (iii) either $P$ is a central product of $m$ copies of
$D_8$ or $P$ is a central
product of $m - 1$ copies of $D_8$ with the quaternion group
$Q_8$.

In both cases $C$ acting faithfully on the elements of $P \backslash Z(P)$ provides 
the required embedding.
\end{proof}

For a finite group $G$ we denote by $R_0(G)$ the smallest degree
of a non-trivial complex projective representation of~$G$.
For simple groups $L$ of Lie type strong lower bounds for
$R_0(L)$ follow from the work of Landazuri and Seitz \cite{LS}
(see \cite[Table 3.3A]{KL} and \cite{Lu}).
See also \cite{KL} and \cite{CCNPW} for the value of $R_0(L)$
when $L$ is an alternating or sporadic group.

Denote by $P(G)$ the minimal degree of a faithful permutation
representation of a group~$G$.
If $L$ is a simple group with a proper subgroup~$H$, then $H$
considered as a subgroup of $\text{\rm Aut}(L)$ is corefree,
hence we have
\[
P\bigl(\text{\rm Aut}(L)\bigr) \leq |L : H |\, |\text{\rm Out}(L)|.
\]

When $L$ is of Lie type and $H$ is a parabolic subgroup of $L$ 
we can improve the above bound as follows: The group $H$ is invariant 
under the subgroups $D, \Phi \leq \mathrm{Aut}(L)$ of diagonal and field 
automorphisms of $L$. 
Then $HD\Phi \leq \mathrm{Aut}(L)$ is a corefree subgroup and hence we have
\[ P\bigl(\text{\rm Aut}(L)\bigr) \leq 6|L : H | \]
if $H$ is parabolic.
Sharp bounds for $P(L)$ when $L$ is a classical simple group can
be found in \cite[Table 5.2A]{KL} and they are achieved for parabolic subgroups.

If $L$ is an exceptional group of Lie type one can easily find a
maximal parabolic subgroup $P$ of small index.
A good lower bound for the order $P$ follows by noting that
$|P|$ is divisible by the order of a Borel subgroup of $G$ and
also by the order of the Levi factor corresponding to $P$
(see \cite[p.~179--181]{KL} for a quick account of these
standard facts).
See also \cite{Wi} for the detailed structure of many of these
maximal subgroups.

Small index subgroups in sporadic simple group can be found in
\cite{CCNPW}.

In the proof of Theorem~\ref{thm:4} we use CFSG via the
following.

\begin{prop}
\label{prop:2.2} There is an absolute constant $c_0$ such that if 
$L$ is a nonabelian finite simple group then
\[
P(\text{\rm Aut}(L)) \leq c_0 R_0(L)^2.
\]
In fact $c_0$ can be taken to be $10^{10}$.
\end{prop}

\begin{proof}
This follows easily from the above mentioned results by inspection.
Note that the bound with $c_0=10^{10}$ is quite sharp for the Monster simple group.
\end{proof}

Recall that a group $H$ is \emph{quasisimple} if it is a perfect,
central extension of a simple group~$L$. It is known \cite{Hu}
that $\text{\rm Aut}(H)$ is isomorphic to a subgroup of $\text{\rm
Aut}(L)$.

The \emph{components} of a group are its subnormal quasisimple subgroups.
The subgroup $E = E(G)$ is generated by the components of~$G$.
The components of $E(G)$ are exactly the components $C_j$ of $G$
and $\text{\rm Aut}(E)$ permutes the components among themselves
\cite{As}.

Denote the orbits of $\text{\rm Aut}(E)$ on the components by
$B_k$ $(k = 1,2, \dots)$.
Then $B_k$ consists of $t_k$ isomorphic components with central
quotient~$L_k$.
The automorphism group of $B_k$ has a natural embedding into\break
$\text{\rm Aut}(L_k)\, wr \, \Sym(t_k)$.
Moreover $\text{\rm Aut}(E)$ has an embedding into $\underset{k}{\Pi} \,
\text{\rm Aut}(B_k)$.
These observations imply the following

\begin{prop}
\label{prop:2.3}
\[
P(\text{\rm Aut}(E)) \leq \sum_j P \bigl(\text{\rm Aut} (C_j/
Z(C_j))\bigr)
\]
where the sum is taken over all components $C_j$ of $G$.
\end{prop}

The generalised Fitting subgroup of $G$ is $F^*(G) = E(G) F(G)$
(where $F(G)$ is the Fitting subgroup of~$G$).
The most significant fact about $F^*(G)$ is that $C_G(F^*(G))
\leq F^*(G)$.
Denote $C_G(F^*(G)) = Z(F^*(G))$ by~$Z$.

The subgroups $E$ and the Sylow subgroups $O_p(G)$ of $F(G)$ are
characteristic in $G$ and their product is $F^*(G)$.
Hence $G/Z$ has a natural embedding into
\[
\text{\rm Aut}(E) \times  \prod\limits_{p/|F(G)|} \text{\rm Aut} (O_p(G)).
\]

Recall that an irreducible linear group $G \leq GL(V)$ is called
 \emph{imprimitive} if the vector space $V$ can be decomposed
into a direct sum $V = V_1 \oplus \dots \oplus V_t$ with $t > 1$, 
such that every element of $G$ permutes the subspaces $V_i$
among themselves, and $G$ is \emph{primitive} if no such
decomposition exists.
By Clifford's theorem any normal subgroup $N$ of a primitive
group is homogeneous, in particular $N$ acts faithfully and
irreducibly on some subspace~$W$ such that $\dim W$ divides
$\dim V$.
For primitive linear groups we prove the following more precise
version of Theorem~\ref{thm:4}.

\begin{thm}
\label{thm:2.4}
Let $G$ be a finite primitive subgroup of $GL(k, \mathbb C)$.
Then $G/Z(G)$ has an embedding into $\Sym(c_0 k^2)$.
\end{thm}

\begin{proof}
It is known \cite{Di} that $Z = Z(G)$ is the unique maximal
normal abelian subgroup of $G$, hence $Z = Z(F^*(G))$.
Moreover,  $Z$ is a group of scalars, hence cyclic.

Let $p$ be a prime such that $O_p(G)$ is not contained in~$Z$.
Then $O_p(G)$ is the product of the Sylow $p$-subgroup $Z_p$ of
$Z$ and an extraspecial $p$-group which is of exponent $p$ in
case $p \neq 2$ \cite[Lemma 1.7]{LMM}.
Therefore the elements of order $p$ (resp. $\leq 4$) in $O_p(G)$
form a characteristic subgroup $R_p$ of $G$ of symplectic type
and $O_p(G) = Z_p \cdot R_p$.

By the remarks preceding the theorem $G/Z$ has a natural
embedding into $\text{\rm Aut}(E) \times \prod\limits_p \text{\rm Aut}(O_p(G))$.
Since conjugation by elements of $G$ stabilizes $R_p$ and fixes
$Z_p$ elementwise we actually have an embedding of $G/Z$ into
$\text{\rm Aut}(E) \times \prod\limits_p C_{\text{\rm Aut}(R_p)}
(Z(R_p))$, where the product is taken over all primes $p$ such that
$O_p(G) \neq Z_p$.

Consider the normal subgroup $N = E \cdot \prod\limits_p R_p$ of~$G$.
The group $N$ may be considered as an irreducible subgroup of $GL(W)$ for
some~$W$ where $\dim(W)$ divides~$k$.
$N$ is a central product of the symplectic type groups $R_p$ and
the components $C_i$ of~$G$ \cite{As}.
Hence there is a decomposition of $W$ into the tensor product of
spaces $\{W_p\}$ and $\{W_i\}$ such that $W_p$ is an irreducible
$R_p$-module for all $p$ and $W_i$ is an irreducible
$C_i$-module for all $i$ \cite[3.7.1 and 3.7.2]{Gor}.

It is clear that $\dim(W_i) \geq R_0(C_i / Z(C_i))$.
We have \[ k \geq \prod\limits_p \dim(W_p) \cdot \prod\limits_i \dim(W_i).\]
Moreover, if $R_p$ has order $p^{2m_p + 1}$ or $p^{2m_p + 2}$,
then $\dim(W_p) = p^{m_p}$ \cite[5.5.5]{Gor}.
Hence we have $k \geq \prod\limits_p p^{m_p} \cdot
\prod\limits_i R_0(C_i / Z(C_i))$ which gives
\[ \sum_p p^{2m_p} + 
\sum_i R_0(C_i / Z(C_i))^2 \leq k^2.\]

But $G/Z$ has a faithful permutation representation of degree at
most
\[
\sum_p P \bigl(C_{\text{\rm Aut}(R_p)} (Z(R_p))\bigr) + P(\text{\rm Aut}(E)) \leq \]
\[ \leq c_0 \bigl( \sum_p p^{2m_p} + \sum_i R_0(C_i/Z(C_i))^2 \bigr) \leq c_0k^2
\]

using Propositions~\ref{prop:2.1}, \ref{prop:2.2} and
\ref{prop:2.3}. Our statement follows.
\end{proof}

\begin{proof}[The end of the proof of  of Theorem~\ref{thm:4}]
If $G \leq GL(k, \mathbb C)$ is irreducible but imprimitive, then
it can be embedded into a wreath product $G_1 \, wr\,  T$ where $G_1$
is a primitive subgroup of $GL(k_0, \mathbb C)$, $T$ is a
transitive subgroup of $\Sym(t)$ and $t k_0 = k$ $[\text{\rm
Sup}]$.
Now $A = Z(G_1)^t$ is an abelian normal subgroup of $G_1 \, wr\,  T$.

We have $\frac{GA}{A} \leq \frac{G_1 \, wr\,  T}{A}$ and by Theorem \ref{thm:2.4} 
$\frac{G_1\, wr\,  T}{A}$ has an embedding \medskip \\ into $\Sym(c_0 k^2_0 t)$.
Hence $\frac{G}{G \cap A}$ has an embedding into $\Sym(c_0k^2)$.
\medskip

Finally, if $G \leq GL(k, \mathbb C)=GL(V)$ (with $V= \mathbb C ^{(k)}$) is an arbitrary finite linear group then $G$ is completely reducible. Hence it embeds into a direct  product $\prod_j GL(V_j)$ where $V= \oplus_j V_j$ is a decomposition of $V$ into irreducible $\mathbb C G$-modules. Let $k_j= \dim_{\mathbb C} V_j$, so that $\sum_j k_j=k$. Our group $G$ acts irreducibly on each $V_j$ and so by the argument above we find a subgroup $A_j \leq G$ such that its image in $GL(V_j)$ is abelian and $G/A_j$ embeds in $\Sym(c_0k_j^2)$. 
Take $A = \cap_j A_j$. It follows that $A$ is abelian and $G/A$ embeds in $\Sym(\sum _j c_0k_j^2) \leq 
\Sym(c_0k^2)$. Theorem \ref{thm:4} follows.
\end{proof}

By a result of Easdown and Praeger \cite{EP} if $G$ is a
subgroup of $\Sym(t)$, then $G/\text{\rm Sol}(G)$ can also be
embedded into $\Sym(t)$ (where $\text{\rm Sol}(G)$ is the
soluble radical of~$G$). 
Hence Theorem~\ref{thm:4} has the
following immediate consequence

\begin{cor}
\label{cor:2.5}
Let $G$ be a finite linear group of degree~$k$.
Then $G/\text{\rm Sol}(G)$ has an embedding into
$\Sym(c_0k^2)$.\hfill $\square$
\end{cor}

In particular $G/\text{\rm Sol}(G)$ is a linear group of degree
at most~$c_0k^2$.
\begin{question} Suppose $G$ is a finite linear group of degree $k$
Is it true that $G/F(G)$ embeds in $\Sym(ck^2)$ for some constant $c$? What about $G/\mathrm{Frat}(G)$?
\end{question}

\section{Additional remarks}
\label{sec:3}

\setcounter{thm}{0}

Following \cite{Gow} let us put some of the results in this
paper into a more general context.

If the minimal degree of a representation of a group $G$ is at
least~$2$, then it is perfect (i.e.\ it is equal to its
commutator subgroup).
Hence it is reasonable to assume that a quasi\-random group is perfect.
For diverse examples of such groups see \cite{HP}.

\begin{thm}
\label{thm:3.1}
Let $G$ be a perfect group of order~$n$.
Then the following statements are polynomially equivalent, in
the sense that if one statement holds for a constant~$c$, then
all others hold with constants that are bounded by a positive
power of~$c$.

\begin{itemize}
\item[\rm (i)] Every representation of~$G$ has degree at
least~$c_1$.

\item[\rm (ii)]
Any product-free subset of $G$ has size at most $\frac{n}{c_2}$.

\item[\rm (iii)]
For any subset $B$ of size at least $\frac{n}{c_3}$ we have $B^3
= G$.

\item[\rm (iv)]
Every proper subgroup of $G$ has index at least~$c_4$.
\end{itemize}
\end{thm}

\begin{proof}
(i) $\implies$ (ii) follows from Proposition~\ref{prop:0} (due
to Gowers).

(i) $\implies$ (iii) follows from Corollary~\ref{cor:1}.

(iv) is an easy consequence of either (ii) or (iii).

(iv) $\implies$ (i) follows from Theorem~\ref{thm:4}.
\end{proof}

Finally, let us point out an application of
Corollary~\ref{cor:1} to permutation groups.
By a deep result of Fulman and Guralnick \cite{FG} if $G$ is a
simple group acting transitively on a set $X$ then the
proportion of fixed point-free permutations in $G$ is at least
$\delta$ for some absolute constant $\delta > 0$.
This implies that if $G$ is large enough, then each element of
$G$ is a product of three fixed point-free permutations.

\section{Appendix: a short proof of a special case of Proposition \ref{prop:0}}
In this section we will give a short version of Gowers' proof of Proposition \ref{prop:0} 
in the case when one of the sets $A,B,C$ is symmetric. This case is enough for most of the applications 
above. We stress that we don't claim originality: certainly 
all the elements of the argument below are already present in Gowers' proof in \cite{Gow}. We believe
that it's worth presenting a simplified version which can fit in one page.
\bigskip

\textbf{Proof:} Suppose $A,B,C$ are three subsets of $G$ one of which coincides with its inverse and 
such that $|A||B||C|>n^3/k$. We have to show that the equation $ab=c$ has a solution with $a \in A, b\in B$ and 
$c \in C$. By cyclically permuting and inverting some of $A,B,C$ we may assume without loss of generality
that $B=B^{-1}$ is the symmetric set.

Let $V= \mathbb C G$ be the group algebra over $\mathbb C$ considered as a complex vector space with 
basis $G$. We will consider $V$ as a left $G$-module equipped the standard Hermitian inner product 
(so the elements $\mathbf{g}$ of $G \subset V$ form an orthonormal basis of $V$). 

Let $X=(x_{g,h})$ be the $n$ by $n$ matrix labelled by $g,h \in G$ such that $x_{g,h}=0$ if 
$h^{-1}g \not \in B$ and $x_{g,h}=1$ otherwise. Then $X$ is a real symmetric matrix defining a linear map 
$X \in \mathrm{End}_{\mathbb C} V$. For each $\mathbf{u} \in G \subset V$ we have $X \mathbf{u}= 
\sum_{b \in B} \mathbf{ub}$ which shows 
that $X$ is in fact a $G$-module endomorphism, i.e. $Xg \mathbf v=gX \mathbf v$ for all $g \in G$ and
$\mathbf v \in V$.

Note that every row and column sum of $X$ is $|B|$, hence $\lambda_1= |B|$ 
is an eigenvalue of $X$ with eigenvector $\mathbf{e}= \sum_{g \in G}\mathbf{g} \in V$. Let $I=\mathbf{e}^{\perp}$ 
be the augmentation ideal of $V$. It is both $X$ and $G$-invariant and clearly 
doesn't have $G$-invariant vectors.

Now, since $X$ is symmetric it has real eigenvalues. Let 
$\lambda \in \mathbb R$ be an eigenvalue of $X$ on $I$ with eigenspace $V_\lambda$. Since $X$ is a $G$-endomorphism 
it follows that $V_\lambda$ is a nontrivial $G$-module and hence $\dim V_\lambda \geq k$.
 
The eigenvalues of $X^2$ are exactly the squares of the eigenvalues of $X$ with 
the same multiplicities. Thus we have that $tr(X^2) \geq k \lambda^2$. But $tr (X^2)=tr (X^t X)$ is exactly 
the sum of all entries of 
$X$ which is $n|B|$. It follows that $\lambda^2 \leq n|B|/k$ holds for all eigenvalues $\lambda$ of $X$ on $I$ and 
therefore 
\begin{equation} \label{eqn} 
|X \mathbf v |^2 \leq \frac{n|B|}{k} |\mathbf v|^2 \textrm{ for all } \mathbf v \in I.
\end{equation}

Let $\mathbf{v}= n \sum_{g \in A} \mathbf{g} \in V$. We can write $\mathbf{v}= \mathbf{v_1}+ \mathbf{v_2}$ where
$\mathbf{v_1}=|A|\mathbf{e}$ and $\mathbf{v_2} \in I$ with $|\mathbf{v_2}|^2=|A|n(n-|A|)<n^2|A|$. Assuming that $ab=c$ has no solution 
$a \in A, b \in B, c \in C$ we deduce from the definition of $X$ that 
$X\mathbf{v} \in \sum_{g \in G \backslash C} \mathbb C \mathbf{g}$.
However 
\[ X \mathbf{v} = X \mathbf{v_1}+ X \mathbf{v_2}= |B||A|\mathbf{e} +X \mathbf{v_2}\]
and it follows that the vector $X \mathbf{v_2}$ has coordinates equal to $-|B||A|$ in at least $|C|$ positions, so
$|X \mathbf{v_2}|^2 \geq |C||A|^2|B|^2$. On the other hand by (\ref{eqn})
\[ |C||A|^2|B|^2 \leq |X \mathbf{v_2}|^2 \leq \frac{n|B|}{k} |\mathbf{v_2}|^2 \leq \frac{n|B|}{k} n^2 |A| \]
which implies $|A||B||C| \leq n^3/k$ and this contradicts the starting assumptions.
Hence $ab=c$ has a solution as stated in the Proposition.

\subsubsection*{Acknowledgement}
We would like to thank Vera  T. S\'os for drawing our attention
to the paper of Gowers, and Aner Shalev for sending us his yet
unpublished paper~\cite{Sh}.


\footnotesize
\begin{thebibliography}{99}

\bibitem[As]{As}
M. Aschbacher,
{\it Finite group theory}, University Press, Cambridge, 1986.

\bibitem[Ba]{Babai}L. Babai,
On the diameter of Eulerian orientations of a graph.
In: {\em Proc. 17th Ann. Symp. on Discrete Algorithms (SODA'06)},
ACM--SIAM 2006, pp. 822-831.

\bibitem[BS]{BS}
L. Babai and V. T. S\'os,
Sidon sets in groups and induced subgraphs of Cayley graphs,
{\it Europ. J. Comb.} {\bf 6} (1985), 101--114.

\bibitem[Co]{Co}
M. Collins,
On Jordan's theorem for complex linear groups, 
\emph{J. Group Theory}, to appear.

\bibitem[CCNPW]{CCNPW}
J. H. Conway, R. T. Curtis, S. P. Norton, R. A. Parker and R. A. Wilson,
{\it An Atlas of finite groups},
University Press, Oxford, 1985.

\bibitem[Di]{Di}
J. D. Dixon,
{\it The structure of linear groups},
Van Nostrand--Reinhold, Princeton, 1971.

\bibitem[Do]{Do}
L. Dornhoff,
{\it Group representation theory. Part A:
Ordinary representation theory},
Pure and Applied Mathematics 7,
Marcel Dekker Inc., New York, 1971.

\bibitem[EP]{EP}
D. Easdown and C. E. Praeger,
On minimal faithful permutation representations of finite groups,
{\it Bull. Aus. M. S.} {\bf 38} (1988), 207--220.

\bibitem[FG]{FG}
J. Fulman and R. M. Guralnick,
Derangements in simple and primitive groups,
in: {\it Groups, Combinatorics and Geometry},
(Durham 2001), World Scientific, New Jersey, 2003, pp. 99--121.


\bibitem[Gor]{Gor}
D. Gorenstein,
{\it Finite groups},
Harper and Row, New York, 1968.

\bibitem[Gow]{Gow}
W. T. Gowers, Quasirandom groups,
preprint.

\bibitem[GR]{GR}
B. Green and I. Z. Ruzsa,
Sum-free sets in abelian groups, preprint.

\bibitem[Gr]{Gr}
R. L. Griess, Jr.,
Automorphisms of extra special groups and nonvanishing degree 2
cohomology,
{\it Pacific J. Math.} {\bf 48} (1973), 403--422.


\bibitem[GL1]{GL1}
R. M. Guralnick and M. Lorenz,
Orders of finite groups of matrices, preprint.

\bibitem[GL2]{GL2}
R. M. Guralnick and F. L\"ubeck,
On $p$-singular elements in Chevalley groups in characteristic
$p$, {\it Groups and Computation. III (Columbus, OH, 1999)},
169--182, Ohio State Univ. Math. Res. Inst. Publ. 8, de Gruyter,
Berlin, 2001.

\bibitem[He]{He}
H. A. Helfgott,
Growth and generation in $SL_2(\mathbb Z / p\mathbb Z)$,
{\it Annals of Math.}, to appear.

\bibitem[HP]{HP}
D. F. Holt and W. Plesken,
{\it Perfect groups},
University Press, Oxford, 1989.

\bibitem[Hu]{Hu}
N. J. S. Hughes,
The structure and order of the group of central automorphisms of
a finite group,
{\it Proc. London Math. Soc.} {\bf 52} (1951), 377--385.

\bibitem[Ke]{Ke}
K. S. Kedlaya,
Product-free subsets of groups,
{\it Amer. Math. Monthly} {\bf 105} (1998), 900--906.

\bibitem[KL]{KL}
P. B. Kleidman and M. W. Liebeck,
{\it The subgroup structure of the finite classical groups},
LMS Lecture Note Series 129,
University Press, Cambridge, 1990.

\bibitem[LS]{LS}
V. Landazuri and G. M. Seitz, On the minimal degrees of projective
representations of the finite Chevalley groups, {\it J. Algebra}
{\bf 32} (1974), 418--443.

\bibitem[La]{La}
M. Larsen,
Word maps have large image,
{\it Israel J. Math.} {\bf 139} (2004), 149--156.

\bibitem[LSh1]{LSh1}
M. Larsen and A. Shalev, Word maps and Waring type
problems, preprint. 

\bibitem[LSh2]{LSh2}
M. W. Liebeck and A. Shalev,
Fuchsian groups, finite simple groups, and representation varieties,
{\it Invent. Math.} {\bf 159} (2005), 317--367.

\bibitem[LSh3]{LSh3} M. W. Liebeck and A. Shalev, Simple groups, probabilistic methods, and a conjecture of Kantor and Lubotzky.  \emph{J. Algebra}  {\bf 184}  (1996),  no. 1, 31--57.

\bibitem[LMM]{LMM}
A. Lucchini, F. Menegazzo and M. Morigi,
On the number of generators and composition length of finite
linear groups,
{\it J. Algebra\/} {\bf 243} (2001), 427--447.

\bibitem[L\"u]{Lu}
F. L\"ubeck,
Smallest degrees of representations of exceptional groups of Lie
type,
{\it Comm. Algebra} {\bf 29} (2001), 2147--2169.

\bibitem[PS]{PS}
A. Penfold Street,
Sum-free sets, in: Springer Lecture Notes in Math. 292,
123--272, 1972.

\bibitem[Sh]{Sh}
A. Shalev,
Word maps, conjugacy classes, and a non-commutative Waring-type
theorem,
{\it Annals of Math.}, to appear.

\bibitem[Su]{Su}
D. A. Suprunenko,
{\it Matrix groups}, A.M.S. Providence, 1976.

\bibitem[Wi]{Wi}
R. A. Wilson,
{\it Finite simple groups}, book in preparation.

\end{thebibliography}
\end{document}